\documentclass{amsart}
\usepackage{amssymb, latexsym, euscript}
\usepackage[usenames,dvipsnames,svgnames,table]{xcolor}

\def\sA{{\mathfrak A}}      
      
   \def\sH{{\mathfrak H}}   
   \def\sK{{\mathfrak K}}   
      
      \def\sR{{\mathfrak R}}

\def\st{{\mathfrak t}}

\def\bB{{\bold B}}

\def\bL{{\bold L}}

\def\bm\chi{\mbox{\boldmath$\chi$}}
\def\half{{\frac{1}{2}}}

\def\ker{{\rm ker\,}}
\def\cker{{\rm \overline{ker}\,}}
\def\ran{{\rm ran\,}}
\def\cran{{\rm \overline{ran}\,}}
\def\dom{{\rm dom\,}}
\def\mul{{\rm mul\,}}

\def\cdom{{\rm \overline{dom}\,}}

\def\dim{{\rm dim\,}}

\let\xker=\ker \def\ker{{\xker\,}}

\DeclareMathOperator{\hoplus}{\, \widehat \oplus \,}

\newtheorem{theorem}{Theorem}[section]
\newtheorem{proposition}[theorem]{Proposition}
\newtheorem{corollary}[theorem]{Corollary}
\newtheorem{lemma}[theorem]{Lemma}

\theoremstyle{definition}
\newtheorem{example}[theorem]{Example}
\newtheorem{remark}[theorem]{Remark}

\numberwithin{equation}{section}

\begin{document}

\title[Non-standard semibounded forms and their representations]{Non-standard behaviour of semibounded forms and their representations}

\author[S.~Hassi]{S.~Hassi}

\address{Department of Mathematics and Statistics \\
University of Vaasa \\
P.O. Box 700, 65101 Vaasa \\
Finland} \email{seppo.hassi@uwasa.fi}

\author[H.S.V.~de~Snoo]{H.S.V.~de~Snoo}

\address{Bernoulli Institute for Mathematics, Computer Science and Artificial Intelligence \\
University of Groningen \\
P.O. Box 407, 9700 AK Groningen \\
Nederland}
\email{h.s.v.de.snoo@rug.nl}

\dedicatory{Dedicated to the memory of our friend Franek Szafraniec (1940-2025)}


\date{\today}


\maketitle

\begin{abstract}
A number of results is presented concerning the extreme behaviour of products of
linear operators or relations from one Hilbert space to another Hilbert space.
It is shown that there exists a relation $T$ from $\sH$ to $\sK$ such that the closure
of the nonnegative relation $T^*T$ in $\sH$ is not equal to the nonnegative selfadjoint relation $T^*T^{**}$;
moreover, this may happen for a singular but also for a regular (i.e., closable) relation $T$.
In fact, as an extreme situation the paper contains a construction of a densely defined closable operator
$T$, such that $\dom T^*T=\{0\}$, while $T^*T^{**}$ is a densely defined selfadjoint operator.
Likewise it is shown that there exist a linear relation $T$ from $\sH$ to
$\sK$ and a nonnegative everywhere defined contraction $R$ in $\sK$
for which the product $RT$ is singular and the product $R^{1/2} T$ is regular.
Such results illustrate situations that can arise in the Lebesgue type decomposition
of linear relations and the Lebesgue type decomposition of semibounded forms.
\end{abstract}

\section{Introduction}

Let $\sH$ and $\sK$ be Hilbert spaces and let $\bL(\sH,\sK)$ the linear space of
all linear relations from $\sH$ to $\sK$ with the usual addition; see \cite{BHS}.
Likewise $\bB(\sH,\sK)$ denotes the linear space of all everywhere defined bounded linear operators.
When $\sK=\sH$ these notations will be abbreviated to $\bf{L}(\sH)$ and $\bf{B}(\sH)$, respectively.
For $T \in \bL(\sH,\sK)$ its adjoint $T^*$ belongs to $\bL(\sK, \sH)$
 and $T^{**} \in \bL(\sH, \sK)$ is the closure of $T$ (as a graph in the Cartesian product $\sH \times \sK$).
The relation $T$  is said to be \textit{regular} (or \textit{closable}) if   $\mul T^{**}=\{0\}$,
so that the closure $T^{**}$ is (the graph of) an operator.
Likewise, the relation $T$ is said to be \textit{singular}
if $T^{**}$ is the Cartesian product of  closed linear subspaces in $\sH$ and $\sK$, respectively.
For a relation $T \in \bL(\sH,\sK)$  there is always a \textit{Lebesgue type decomposition}
as a sum:
\begin{equation}\label{leb1}
 T=T_1+T_2,
\end{equation}
where $T_1 \in \bL(\sH,\sK)$ is a regular (closable) operator and $T_2 \in \bL(\sH,\sK)$
is a singular relation.  In particular, these summands are obtained via
\begin{equation}\label{leb2}
 T_1=(I-K) T \quad \mbox{and} \quad T_2=KT,
\end{equation}
where $K \in \bB(\sK)$ is a nonnegative contraction, so that $T_1$ and $T_2$ have
the desired properties; see \cite{HS2023a}.  The components are mutually singular
if and only if $K\in \bB(\sK)$ is an orthogonal projection,
 so that $\ran T_1 \perp \ran T_2$. The particular case of \eqref{leb1} and \eqref{leb2}
where $K$ is the orthogonal projection onto
the multivalued part $\mul T^{**}$ leads to the standard \textit{Lebesgue decomposition},
which is usually denoted by
\begin{equation*}
 T=T_{\rm reg}+T_{\rm sing}.
\end{equation*}
Here $T_{\rm reg}$ is in fact the largest closable operator that is dominated by $T$,
see \cite{HSSS2007}. In fact, the regular part $T_{\rm reg}$ can be written as
\begin{equation}\label{regu}
 T_{\rm reg}=(I-P) T,
\end{equation}
where $P$ is the orthogonal projection in $\sK$ onto $\mul T^{**}$.

\medskip

Semibounded sesquilinear forms $\st$ on a Hilbert space $\sH$
give rise to decompositions similar to \eqref{leb1}
via an associated representing operator $T \in \bL(\sH, \sK)$
with some Hilbert space $\sK$ via
\[
 \st[f,g]=(Tf,Tg), \quad f,g \in \dom \st=\dom T.
\]
The form is not assumed to be densely defined in which case $T^*T$ is
a nonnegative linear relation with a nontrivial multivalued part.
Moreover, the form $\st$ is closed (closable or singular)
if and only if the operator $T$ is closed (closable or singular, respectively); see \cite[Lemma~2.3]{HS2023b}.
In particular, if the form $\st$ is closed then $T^*T$ is a selfadjoint relation.
For the decomposition of a semibounded form
one now has to consider a nonnegative contraction $K \in \bB(\sK)$ for which $(I-K)^\half T$
is regular (closable) and $K^\half T$ is singular; cf.  \cite{HS2023b}.
Of course, when $K$ is an orthogonal projection,
this leads to the same description as for \eqref{leb2}.
However, in general, the relations $(I-K)T$ and $(I-K)^\half T$
need not be regular simultaneously and, likewise, $KT$ and $K^\half T$
need not be singular simultaneously.
Closely connected with the above forms are the products $T^*T$ and $T^*T^{**}$
when $T \in \bL(\sH, \sK)$ is a linear relation. It is clear that the relation
$T^*T$ is nonnegative, whereas the relation $T^*T^{**}$ is selfadjoint.
The inclusion
\begin{equation*}\label{leb4}
 T^*T  \subset T^*T^{**}
\end{equation*}
is straightforward. The relation $T^*T$ need not be closed
and its closure need not be selfadjoint.
Note that $(T^*T)^{**}$ is selfadjoint
if and only if  it is equal to $T^*T^{**}$.

\medskip

This paper will illustrate some of the above phenomena
by rather concrete examples.  There have been quite a few papers where related results have
been given, see \cite{DM23}, \cite{GS0}, \cite{GS1}, \cite{Mo0},  \cite{Mo2},
and \cite{ST12, ST14, ST16, ST19, ST21}, and there is even a book
with some counterexamples to desired statements \cite{Mo1}.
Here the approach uses a construction that goes back to Von Neumann
\cite{vN}  and which has been refined in \cite{FW}.
It involves a pair of nonnegative injective bounded linear operators in a
Hilbert space, whose ranges have a trivial intersection, whereas the ranges
of their square roots coincide. Related useful ideas can of course be found in \cite{Doug}.
The work on Lebesgue type decompositions has given rise to some of the present results,
\cite{An}, \cite{J}, \cite{Si}.
In the context of linear relations the study of Lebesgue type decompositions goes back about twenty years
with some joint work \cite{HSSS2007} with Zolt\`an Sebesty\'en and Franek Szafraniec,
which was continued with Franek in \cite{HSS}.

\section{Some preliminaries}\label{prelim}

Let $A, B \in \bB(\sH,\sK)$, where $\sH$ and $\sK$ are Hilbert spaces.
In this section the inclusion $\ran A \subset \ran B$ will be reviewed
in terms of the Douglas lemma in \cite{Doug} and
an important construction of Von Neumann from \cite{FW} will be recalled.

\begin{lemma}\label{Douglas}
 Let $A, B \in \bB(\sH,\sK)$.
 Then the following statements are equivalent:
\begin{enumerate}\def\labelenumi{\rm (\roman{enumi})}
\item $\ran A \subset \ran B$;
\item $A=BW$ for some  $W \in \bB(\sH)$;
\item $AA^{*} \le \lambda BB^{*}$ for some $\lambda \ge 0$.
 \end{enumerate}
The operator $W$ can be chosen such that $\ran W \subset (\ker B)^\perp$, in which case
its restriction to $(\ker A)^\perp$ is unique.
 \end{lemma}

It is a direct  consequence of Lemma \ref{Douglas} that the following statements are equivalent:
\begin{enumerate}\def\labelenumi{\rm (\roman{enumi})}
\item [(i')] $\ran A = \ran B$;
 \item [(ii')] $AA^{*} \le \lambda BB^{*} \leq \rho AA^*$ for some $\lambda, \rho \ge 0$.
\end{enumerate}

\medskip

Observe that any $W \in \bB(\sH)$ in (ii) of Lemma \ref{Douglas}
has the property that $W$ maps $\ker A$ boundedly into $\ker B$
and that this part of $A$ may be completely arbitrary. In the context of
Lemma \ref{Douglas} one may of course define $Wf=0$ for all $f \in \ker A$,
in which case
\[
\|W\|^2 =\inf\,\{\mu \geq 0 :\,AA^*\leq\mu BB^*\}.
\]
However, in the following corollary this flexibility in the choice of $W$ is used in a different way.
A proof is included for convenience; for some closely related statements, see \cite[App.~D]{BHS}.

\begin{corollary}\label{Douglas1}
 Let $A, B \in \bB(\sH,\sK)$.
Then the following statements are equivalent:
\begin{enumerate}\def\labelenumi{\rm (\roman{enumi})}
\item $\ran A = \ran B$ and $\dim \ker A=\dim \ker B$;
\item $A=BW$ for some $W \in \bB(\sH)$ with $W^{-1} \in \bB(\sH)$.
\end{enumerate}
In particular, if $A,B\in \bB(\sH)$ are selfadjoint then {\rm (i)} simplifies to $\ran A = \ran B$.
 \end{corollary}

\begin{proof}
It is clear from Lemma \ref{Douglas} that (ii) implies (i).
Conversely, if (i) holds, then one can choose bounded linear operators $W, Z$ such that
$A=BW$ and $B=AZ$ where $W$ maps $(\ker A)^\perp$ into $(\ker B)^\perp$,
while $Z$ maps $(\ker B)^\perp$ into $(\ker A)^\perp$.
Observe that $ZW$ is the identity on $(\ker A)^\perp$
while $WZ$ is the identity on $(\ker B)^\perp$. Thus $W$ maps $(\ker A)^\perp$ bijectively onto $(\ker B)^\perp$.
It remains to extend the definition of $W$ to all of $\sH$. Since $\dim \ker A=\dim \ker B$ one may choose $W$ to
be a bounded and boundely invertible operator from $\ker A$ onto $\ker B$. This shows that (ii) holds.
\end{proof}

\medskip

It is not difficult to see that (i), (ii), and (iii) in Lemma \ref{Douglas} are equivalent to
\[
 AA^*=BMB^* \quad \mbox{for some nonnegative} \quad M \in \bB(\sH).
\]
The operator $M$ may be chosen such that $\ran M \subset (\ker B)^\perp$,
in which case the restriction of $M$ to  $(\ker B)^\perp$  is unique. Moreover,
the conditions $\ran A = \ran B$ and $\dim \ker A=\dim \ker B$ are equivalent to
$M \in \bB(\sH)$ and $\ran M=\sH$ (or, alternatively,  $M^{-1} \in \bB(\sH)$).

\medskip

Here are some consequences for nonnegative operators in $\bB(\sH)$.

\begin{lemma}\label{surprise}
Let $A, B \in \bB(\sH)$ be nonnegative.
Then the following statements are equivalent:
\begin{enumerate}\def\labelenumi{\rm (\roman{enumi})}
\item $\ran A=\ran A^{1/2}$;
\item $\ran A$ is closed.
\end{enumerate}
Moreover, the following statements are equivalent:
\begin{enumerate}\def\labelenumi{\rm (\roman{enumi})}\setcounter{enumi}{2}
\item $\ran A^{\half}=\ran B^{\half}$;
\item $A=B^{\half}MB^{\half}$ with $M \in \bB(\sK)$ nonnegative and $\ran M=\sK$.
\end{enumerate}
\end{lemma}

The following well-known result goes back to Von Neumann; see \cite[p.~273]{FW}.
A proof due to Dixmier can be found in \cite[Theorem~3.7]{FW}.

\begin{lemma}\label{fill0}
There exist nonnegative operators
$A,B \in \bB(\sH)$, such that
\[
\ran A\cap\ran B=\{0\}.
\]
\end{lemma}

The result in Lemma \ref{fill0} has been modified by Fillmore and Williams; see \cite[p.~278]{FW}.
Their simple proof, using Lemma \ref{Douglas}, has been included.

\begin{lemma}\label{fill}
There exist nonnegative operators
$A,B \in \bB(\sH)$, such that
\[
\left\{
\begin{array}{l}
\ran A\cap\ran B=\{0\}, \\
\ran A^{\half}=\ran B^{\half}.
\end{array}
\right.
\]
Consequently, without loss of generality one may assume that $A$ and $B$ are injective.
\end{lemma}

\begin{proof}
Let the operators $A_1$ and $A_2$ in $\bB(\sH)$
be nonnegative such that they satisfy
$\ran A_1 \cap \ran A_2=\{0\}$.
Then introduce $A$ and $B$ in $\bB(\sH)$ by
 \[
 A=A_1+A_2 \quad \mbox{and} \quad B=A_1+2A_2.
\]
By the assumption one sees that
$\ran A \cap \ran B=\{0\}$.
Furthermore, the definition of $A$ and $B$ implies
\[
 A \leq B \leq 2A.
\]
These inequalities imply as a direct consequence of Lemma \ref{Douglas} that
 there is the identity $\ran A^\half=\ran B^\half$. Furthermore, one has $\ker A=\ker B$
 and without loss of generality it may be therefore be assumed that $\ker A=\ker B=\{0\}$.
\end{proof}

\section{Closable and singular relations}

Let $T \in \bL(\sH,\sK)$, where $\sH$ and $\sK$ are Hilbert spaces.
It has been remarked that $T$  is said to be regular (or closable) if   $\mul T^{**}=\{0\}$
and that $T$ singular if $T^{**}$ is the Cartesian product of  closed linear subspaces in $\sH$ and $\sK$.
Hence $T$ is singular if and only if  $\dom T \subset \ker T^{**}$ or,
equivalently, $\ran T \subset \mul T^{**}$. It is also clear that $T$ is singular if and only if $T^*$
is singular, in other words if $\dom T^* \subset \ker T^*$ or, equivalently, $\ran T^* \subset \mul T^*$.
Recall the Lebesgue type decomposition of $T \in \bL(\sH,\sK)$ in \eqref{leb1} and \eqref{leb2}.
This section will contain some straightforward statements about the products in \eqref{leb2}.

\medskip

Let $T \in \bL(\sH,\sK)$  and let $R \in \bB(\sK)$.
The interest is in some properties of the product $RT$:
\[
 RT=\big\{ \{f, Rf'\} :\, \{f, f'\} \in T \big\}.
\]
Note that $RT \in \bL(\sH, \sK)$ and
\[
 \dom RT=\dom T \quad \mbox{and} \quad \mul RT=R\, \mul T.
\]
Recall from \cite[Proposition 1.3.9]{BHS} that $(RT)^*=T^* R^*$, so that
\begin{equation}\label{Aa}
 RT^{**} \subset (RT)^{**} \quad \mbox{and} \quad R\, \mul T^{**} \subset \mul (RT)^{**}.
\end{equation}

\begin{lemma}\label{simple}
Let $T \in \bL(\sH, \sK)$ be a linear relation
and let $R \in\mathbf{B}(\sK)$.
Then the following statements hold:
\begin{enumerate}\def\labelenumi{\rm(\alph{enumi})}
\item if $RT$ is closable and $\ker R \cap \mul T^{**}=\{0\}$, then $T$ is closable;
\item if $T$ is singular, then $RT$ singular.
\end{enumerate}
\end{lemma}

\begin{proof}
(i) Assume $RT$ is closable, i.e., $\mul (RT)^{**}=\{0\}$. Then \eqref{Aa} shows that
$\mul T^{**} \subset \ker R$. With the assumption $\ker R \cap \mul T^{**}=\{0\}$
this gives $\mul T^{**}=\{0\}$. Thus, $T$ is closable.

(ii) Assume that $T$ is singular or, equivalently, that $T^*$ is singular, in other words
$\dom T^* \subset \ker T^*$. This clearly implies $\dom T^*R^* \subset \ker T^*R^*$.
Now observe that the identity $T^*R^*=(RT)^*$ implies that $RT$ is singular.
\end{proof}

A more symmetric statement can be found in the next corollary.

\begin{corollary}\label{nieuw}
Let $T \in \bL(\sH,\sK)$ be a linear relation
and let $R\in\mathbf{B}(\sK)$ be nonnegative.
Then the following statements hold:
\begin{enumerate} \def\labelenumi{\rm(\alph{enumi})}
\item if $RT$ is closable, then $R^{\half} T$ is closable;
\item if $R^\half T$ is singular, then $RT$ is singular.
\end{enumerate}
\end{corollary}

\begin{proof}
(a) Assume that $RT=R^\half (R^\half T)$ is closable and observe that
\[
 \ker R^\half \cap \mul (R^\half T)^{**}
 =\ker R^\half \cap R^\half (\mul T^{**})\subset \ker R^{1/2} \cap \ran R^{1/2}=\{0\}.
\]
 Thus Lemma \ref{simple} gives (a).

 (b) This assertion is an immediate consequence of Lemma \ref{simple}.
\end{proof}

The implications in (a) and (b) in Corollary \ref{nieuw} cannot be reversed as will be
shown in the next section.

\section{Some constructions of regular and singular relations}

The representation of linear relations in terms of pairs of bounded linear operators
can be found in \cite{HSS2018}, where also some examples were given.
In this section it will be shown that nonnegative operators $A, B \in \bB(\sH)$ for which
$\ran A \cap \ran B=\{0\}$, as given in Section \ref{prelim},
can be used to illustrate some of the above mentioned issues in the theory of
closable and singular operators or relations.
First it will be pointed out that such operators provide
a very simple, but extreme, example of a singular operator; cf. Lemma \ref{fill}.
For examples of a different nature, see \cite{HS2022}.

\begin{proposition}\label{einz}
Let $A, B \in \bB(\sH)$  be injective nonnegative operators.
Then the linear relation $T \in \bL(\sH)$, given by
\begin{equation}\label{een0}
  T=BA^{-1},
\end{equation}
is a densely defined injective operator.
Moreover, the adjoint $T^*$ is a well-defined closed injective operator given by
\begin{equation}\label{twee}
 T^*=A^{-1}B.
\end{equation}
 Consequently, if in addition $\ran A \cap \ran B=\{0\}$, then $T$
is singular.  In fact, $T^*=\{0\} \times \{0\}$ and $T^{**}=\sH \times \sH$,
in particular, the operator $T$ is not closed or closable.
\end{proposition}

\begin{proof}
By the definition of $T=BA^{-1}$ in \eqref{een0} the identity
\begin{equation}\label{een}
  T=BA^{-1}=\big\{ \{A \varphi, B \varphi\} :\, \varphi \in \sH \big\}
\end{equation}
is clear.
It follows from \eqref{een} that $\dom T=\ran A$, so that $T$ is densely defined.
Moreover, if $A\varphi=\{0\}$, then $\varphi=0$, so that $B\varphi=0$.
Hence, $T$ is the graph of an operator.
It is immediately clear from \eqref{een} that  $T^*$ satisfies
\begin{equation}\label{twee1}
 T^*=A^{-1}B=\big\{ \{f,g\} \in \sH \times \sH :\, Bf=Ag \big\},
\end{equation}
cf. \cite[Proposition 1.3.9]{BHS}.
Observe, that \eqref{twee1} implies
\begin{equation}\label{drie}
\dom T^*=\{ f \in \sH :\, Bf \in \ran A\} \quad \mbox{and}  \quad \ran T^*=\{ g \in \sH :\, Ag \in \ran B\}.
\end{equation}
Thus, if  $\ran A \cap \ran B=\{0\}$,
then \eqref{drie} shows that $T^{*}=\{0\} \times \{0\}$ or $T^{**}=\sH \times \sH$.
Therefore, $T \in \bL(\sH)$ is singular.
\end{proof}

In particular, the conclusions of Proposition \ref{einz} hold when the pair of
injective operators $A, B \in \bB(\sH)$ satisfy the conditions of Lemma \ref{fill}.

\medskip

Now recall the implications stated in (a) and (b) of Corollary \ref{nieuw}.
The next proposition shows that these implications cannot be reversed.

\begin{proposition}\label{mmm}
Let $A, B \in \bB(\sH)$  be nonnegative injective operators with
\begin{equation}\label{ab}
 \ran A \cap \ran B=\{0\} \quad \mbox{and} \quad \ran A^\half=\ran B^\half.
\end{equation}
Then the densely defined operator $T \in \bL(\sH)$ and the operator $R \in \bB(\sH)$, given by
\begin{equation}\label{ab1}
 T=A^{-1} \quad \mbox{and} \quad R=B
\end{equation}
are both nonnegative, and they have the following properties
 \begin{equation}\label{ab2}
 R^{\half}T \quad \textrm{is closable \; and} \quad  RT \quad \textrm{is singular.}
\end{equation}
 \end{proposition}

\begin{proof}
Observe that the operator $RT$ is given by $RT= BA^{-1}$.
By Proposition~\ref{einz} the operator $BA^{-1}$ is singular, which is the second assertion in \eqref{ab2}.
Next note that the operator $R^{\half}T$ can be expressed by
$R^\half T=B^\half A^{-1}$, where $B^\half \in \bB(\sH)$, so that
\[
 (R^\half T)^*=A^{-1} B^\half.
\]
Therefore, one obtains by means of Lemma \ref{surprise}
\[
\begin{split}
&\dom (R^{\half}T)^* =\dom A^{-1}B^{\half} =\ran B^{-\half}A \\
& \hspace{2cm}=\ran B^{-\half}B^{\half}MB^{\half}=\ran MB^{\half}.
\end{split}
\]
This linear subspace is dense in $\sH$, since $\cran B^{\half}=\sH$ and $\ker M=\{0\}$.
Thus $(R^{\half}T)^*$ is densely defined and, hence, the operator $R^{\half} T$ is closable,
which is the first assertion in \eqref{ab2}.
\end{proof}

Of course one may assume without loss of generality that $B$ in \eqref{ab} is
a contraction. Hence there exists a pair in \eqref{ab1}, consisting of a nonnegative operator $T$
and a nonnegative contraction $R$, for which the assertions in \eqref{ab2} hold.

\section{On the product $T^*T$ for a linear relation $T$}

If $T \in \bL(\sH, \sK)$ is a linear relation,
then the product relation $T^*T \in \bL(\sH)$
is nonnegative with
\[
 \mul T^*T=\mul T^*=(\dom T)^\perp
\]
and, moreover, the nonnegative relation $T^*T^{**}$
is selfadjoint; see e.g. \cite[Proposition~1.3.2, Lemma~1.5.8]{BHS}.
Recall that $T^*T^{**}$ decomposes as a componentwise
orthogonal sum
\[
 T^*T^{**}=(T^*T^{**})_{\rm reg} \hoplus (\{0\} \times \mul T^*),
\]
where the orthogonal operator part $(T^*T^{**})_{\rm reg}$ is
a nonnegative selfadjoint operator in $\sH \ominus \mul T^*$;
cf. \cite[Theorem~1.5.1]{BHS}. For the definition of the
regular part see \eqref{regu}.
 Moreover, there is the following useful identity
\begin{equation}\label{regreg}
 T^*T^{**}=(T_{\rm reg})^*(T_{\rm reg})^{**},
\end{equation}
cf. \cite{HS2023seq}.

\medskip

In general, it follows from $T \subset T^{**}$ that the nonnegative relation
$T^*T$ has a nonnegative selfadjoint extension $T^* T^{**}$.
Hence, it is clear that
\begin{equation}\label{teetee2}
T^*T   \subset T^*T^{**} \subset (T^*T)^*.
\end{equation}
In particular, the inclusion \eqref{teetee2} shows that the closure of $T^*T$ satisfies
\begin{equation}\label{teetee3}
(T^*T)^{**}\subset T^*T^{**}.
\end{equation}
In general, there is no equality in \eqref{teetee3}.

\medskip

First it is  shown when there is equality in \eqref{teetee3}
for the case that $T \in \bL(\sH, \sK)$ is a singular relation.

\begin{lemma}\label{Sing}
Let $T \in \bL(\sH, \sK)$ be a singular relation. Then
\begin{enumerate} \def\labelenumi{\rm(\alph{enumi})}
\item $T^*T=\ker T \times (\dom T)^\perp$;
 \item  $T^*T^{**}=\cdom T \times (\dom T)^\perp$.
\end{enumerate}
Consequently, there are the following equivalences:
\begin{enumerate} \def\labelenumi{\rm(\alph{enumi})}\setcounter{enumi}{2}
\item $T^*T=T^*T^{**}$ if and only if $\ker T=\cdom T$;
\item $(T^*T)^{**}=T^*T^{**}$ if and only if $\cker T=\cdom T$.
\end{enumerate}
 \end{lemma}

\begin{proof}
(a)  Since the relation $T$ is singular, one has
$\dom T^*=\ker T^*$.
Now let $\{f,f'\} \in T^*T$, so that  $\{f,h\} \in T$ and $\{h,f'\} \in T^*$
for some element $h \in \sK$.  In fact,
\[
h \in \ran T \cap \dom T^* =\ran T \cap \ker T^*=\{0\}.
\]
Thus $f \in \ker T$ and $f'\in \mul T^*$, which gives $\{f,f'\} \in \ker T \times \mul T^*$.
Therefore $T^*T \subset \ker T \times \mul T^*$. The reverse inclusion is obvious.
Thus (a) has been shown.

(b) Since the relation $T$ is singular, one has $T=T_{\rm sing}$ or, equivalently,
$T_{\rm reg}$ is the zero operator. Hence, $T_{\rm reg}=\dom T \times \{0\}$, which gives
\[
(T_{\rm reg})^{*}=\sK \times (\dom T)^\perp \quad \mbox{and} \quad
(T_{\rm reg})^{**}= \cdom T \times \{0\},
\]
and thus it follows that
\[
(T_{\rm reg})^*(T_{\rm reg})^{**}=\cdom T \times (\dom T)^\perp.
\]
Now apply \eqref{regreg} to obtain (b).
\end{proof}

Lemma \ref{Sing} gives a nice characterization of singularity as a byproduct.

\begin{corollary}
Let $T \in \bL(\sH, \sK)$ be a linear relation. Then the following statements are equivalent:
\begin{enumerate} \def\labelenumi{\rm(\roman{enumi})}
\item $T$ is singular;
\item $T^*T^{**}$ is singular.
\end{enumerate}
In this case $T^*T^{**}=\sA \times \sA^\perp$ for some closed linear subspace $\sA \subset \sH$.
\end{corollary}

\begin{proof}
First make the trivial observation
that the selfadjoint relation $T^*T^{**}$ is singular if and only if $T^*T^{**}=\sA \times \sA^\perp$
for some closed linear subspace $\sA \subset \sH$.

(i) $\Rightarrow$ (ii)  Assume that $T$ is singular. Then it follows from Lemma \ref{Sing} that
$T^*T^{**}=\cdom T \times (\dom T)^\perp$, and, in particular, $T^*T^{**}$ is singular.

(ii) $\Rightarrow$ (i) Assume that $T^*T^{**}$ is singular, i.e., $T^*T^{**}=\sA \times \sA^\perp$
for some closed linear subspace $\sA \subset \sH$.
 Therefore, one obtains
\[
\left\{
\begin{array}{l}
 \mul T^*=\mul T^*T^{**}=\sA^\perp, \\
 \ker T^{**}=\ker T^*T^{**}=\sA.
\end{array}
\right.
\]
Thus $\cran T^*= (\ker T^{**})^\perp= \mul T^*$, i.e., $T$ is singular.
\end{proof}

Now it is not difficult to show that Proposition \ref{einz} leads to
a densely defined singular operator for which the inclusion \eqref{teetee3} is strict.

\begin{proposition}\label{5.3}
Let $A, B \in \bB(\sH)$  be nonnegative injective operators such that
\[
 \ran A \cap \ran B=\{0\}
\]
and let $T=BA^{-1}$ as in Proposition {\rm \ref{einz}}.
Then $T$ is a densely defined injective singular operator and
\begin{equation}\label{ab02}
 T^*T=\{0\} \times \{0\} \quad \mbox{and} \quad T^*T^{**}=\sH \times \{0\},
\end{equation}
so that $T^*T$ and $T^*T^{**}$ are both singular and regular.
In particular,  the inclusion $(T^*T)^{**} \subset T^*T^{**}$ is strict.
\end{proposition}

\begin{proof}
According to Proposition \ref{einz} $T=BA^{-1}$ is a densely defined injective operator
which is singular and satisfies $T^*=\{0\} \times \{0\}$ and $T^{**}=\sH \times \sH$.
This leads to \eqref{ab02} by direct calculation of the products.
It is clear from \eqref{ab02} that the inclusion $(T^*T)^{**} \subset T^*T^{**}$ is strict;
in fact, this also follows from Lemma \ref{Sing}.
\end{proof}

Next an essentially stronger result is proved.
Namely, the inclusion in \eqref{teetee3} may be strict:
there exists a closable operator $T$ for which $\dom T^*T=\{0\}$.

\begin{theorem}\label{thm1}
 Let $A, B \in \bB(\sH)$ be nonnegative injective operators such that
\[
 \ran A \cap \ran B=\{0\} \quad \mbox{and} \quad \ran A^\half=\ran B^\half,
\]
and let $T \in \bL(\sH)$ be given by  $T=B^{\half}A^{-1}$.
Then $T$ is a densely defined closable operator, such that
\begin{equation}\label{T1}
T^*T=\{0\} \times \{0\} \quad \mbox{and} \quad  T^*T^{**}=A^{-1}B^{\half} (A^{-1}B^{\half})^*.
\end{equation}
In particular, the inclusion $(T^*T)^{**} \subset T^*T^{**}$ is strict.
\end{theorem}

\begin{proof}
Recall that $\dom T=\ran A$, so that $T$ is densely defined.
Moreover, since $B^{\half}\in\mathbf{B}(\sH)$
the adjoint of $T$ is an injective operator given by
\[
 T^*=A^{-1}B^{\half}.
\]
Next observe that $T^*=B^{-\half}M^{-1}$ with $M^{-1} \in \bB(\sH)$, see Lemma \ref{surprise}, so that
\[
 \dom T^* =\ran MB^{\half}.
\]
This is a dense set in $\sH$, since $\ran B^{\half}$ is dense and
$M \in \bB(\sH)$ has dense (even closed) range.
Therefore, $T^*$ is densely defined.
Thus $T$ is a closable operator.
Moreover, it is clear that
\[
T^*T=A^{-1}B^\half B^\half A^{-1}= A^{-1}BA^{-1}.
\]
On the other hand, since  $\ran A\cap\ran B=\{0\}$, one has $\dom A^{-1}B=\{0\}$ and,
consequently, $\dom T^*T=\{0\}$ follows.
Observe that $T$ is densely defined so that $T^*T$ is an operator;
hence $T^*T=\{0\} \times \{0\}$.
\end{proof}

Observe that the relation $T \in \bL(\sH, \sK)$ is densely defined if and only if the selfadjoint relation
$T^*T^{**}$ is an operator. This is independent of $T$ being closable.  Therefore, in particular,
if $T \in \bL(\sH,\sK)$ is a densely defined linear operator, then $T^*T$ is an operator that has a selfadjoint
extension $T^*T^{**}$, see \eqref{teetee2}, which is automatically an operator;  cf. \cite{ST21}.

\medskip

The inverse of the operator in Proposition \ref{thm1}
gives also following interesting result.

\begin{proposition}\label{thm11}
Let $A, B \in \bB(\sH)$  be nonnegative injective operators such that
\[
 \ran A \cap \ran B=\{0\} \quad \mbox{and} \quad \ran A^\half=\ran B^\half,
\]
and let $T \in \bL(\sH)$ be given by  $T=AB^{-\half}$.
Then $T$ is a densely defined operator, which is  bounded {\rm(}and hence closable{\rm )},
such that
\[
 T^{**} \in\mathbf{B}(\sH) \quad \mbox{and} \quad \dom TT^*=\{0\}.
 \]
\end{proposition}

\begin{proof}
By definition $\dom T=\ran B^\half$, so that $T$ is densely defined.
Moreover, its adjoint is given by $T^* =B^{-\half}A$, see \cite[Proposition 1.3.9]{BHS}.
Observe that  $\ran A\subset \ran A^{\half}=\ran B^{\half}$ implies that  $\dom T^*=\sH$.
By the closed graph theorem $T^*\in \mathbf{B}(\sH)$,
and therefore also $T^{**}=(B^{-\half}A)^*\in\mathbf{B}(\sH)$.
However,
\[
 TT^*=AB^{-\half}B^{-\half}A=AB^{-1}A,
\]
and $\dom TT^*=\{0\}$, since $\dom B^{-1}A=\{0\}$ due to $\ran A\cap \ran B=\{0\}$.
\end{proof}

\begin{remark}
Let $T$ be the operator from Theorem~\ref{thm1}. Then
\[
 \st[f,g]=(Tf,Tg), \quad f,g \in \dom \st=\dom T,
\]
defines a densely defined closable form and in the analog of Kato's representation theorem \cite{Kato}
one is dealing with the symmetric operator $T^*T$. In fact, the selfadjoint operator $T^*T^{**}$ is also
associated with the form $\st$ in the following extended sense:
 the graph of $T^*T^{**}$
consists of all  $\{f,g\} \in \sH^2$ for which there exists a sequence
$f_n \to_\st f$ such that
\[
 \lim_{n \to \infty}  \st[f_n,h] = (g,h),
\]
for all $h \in \dom \st=\dom T$. For further results along these lines,
see \cite{HS2023b}; for the context of sectorial forms, see \cite{AE}
and \cite{HSW}.
\end{remark}

\section{On the products $T^*T$ and $TT^*$ for a linear relation $T$}

The block relations or block operators  treated in this section often appear in the literature
to produce symmetric and selfadjoint relations or operators, and for instance,
to connect spectral properties of $T^*T$ and $TT^*$, see e.g. \cite{HM2001}.
Such results are also of interest for Pontryagin and Kre\u{\i}n space operators,
like simultaneous $J$-contractivity of $T$ and $T^*$.
Block relations in the present context are used for more flexibility; cf. \cite{HLS20}.

\medskip

Let $T \in \bL(\sH, \sK)$ be a linear relation and
define the linear relation $Z \in \bL(\sK \oplus \sH)$ by
the block relation
\begin{equation}\label{A}
 Z=\begin{pmatrix}
     0 & T \\
     T^* & 0
   \end{pmatrix}.
\end{equation}
For block relations and their calculus, see \cite{HLS20}.
It is easy to check that
\begin{equation}\label{Astar}
 Z^*=\begin{pmatrix}
     0 & T^{**} \\
     T^* & 0
   \end{pmatrix}
\quad \mbox{and} \quad
 Z^{**}=\begin{pmatrix}
     0 & T^{**} \\
     T^{*} & 0
   \end{pmatrix},
\end{equation}
so that $Z$ is essentially selfadjoint.
 By \eqref{A} and \eqref{Astar} one sees that
\begin{equation}\label{A1}
Z^*Z
   =\begin{pmatrix}
     T^{**}T^* & 0 \\
     0 & T^*T
   \end{pmatrix}
   \quad
   \mbox{and}
   \quad
   ZZ^*=\begin{pmatrix}
     TT^* & 0 \\
     0 & T^*T^{**}
   \end{pmatrix}
\end{equation}
are nonnegative relations,  whose domains are given by
\[
\left\{
\begin{array}{l}
\dom Z^*Z=\dom T^{**}T^{*}  \oplus \dom T^*T, \\
\dom ZZ^*=\dom TT^* \oplus \dom T^*T^{**}.
\end{array}
\right.
\]

\medskip

The above construction can  be used to prove the
next result which goes back to Sebesty\'en and Tarcsay \cite{ST14}.
Gesztesy and Schm\"udgen \cite{GS0, GS1, Sch} were the first to give this proof.

\begin{proposition}\label{oneone}
Let $T \in \bL(\sH, \sK)$ be a linear relation.
Then the following statements are equivalent:
\begin{enumerate} \def\labelenumi{\rm(\roman{enumi})}
\item $T$ is closed;
\item $T^*T$ and $TT^*$ are selfadjoint.
\end{enumerate}
\end{proposition}

\begin{proof}
(i) $\Rightarrow$ (ii) This is clear.

(ii) $\Rightarrow$ (i) Note that by assumption the relation
\[
 Z^2=\begin{pmatrix}
    TT^* & 0 \\
     0 & T^*T
  \end{pmatrix}
\]
is selfadjoint.  Since $T^*T$ and $TT^*$ are assumed to be selfadjoint, one has that
$\ran (TT^*+I_\sK)=\sK$ and $\ran (T^*T+I_\sH)=\sH$. Hence
$\ran (Z^2+I_\sK \oplus I_\sH)=\sK \oplus \sH$, which implies that
\[
 \ran (Z \pm i I_\sK \oplus I_\sH)=\sK \oplus \sH.
\]
Since $Z$ is symmetric, it follows that $Z$ is selfadjoint and therefore $T$ is closed.
\end{proof}

The absence of closedness for the linear relation $T$
may lead to extreme situations for $T^*T$ and $TT^*$;  cf. Proposition \ref{oneone}.
Recall for this purpose  Proposition~\ref{5.3}.

\begin{proposition}
The densely defined linear operator $T \in \bL(\sH)$ in Proposition  {\rm \ref{5.3}} is not closed and
there are the following identities
\[
 T^*T=\{0\} \times \{0\}   \quad \mbox{and} \quad TT^*=\{0\} \times \{0\}.
\]
However, $T^*T^{**}$ and $T^{**}T^*$ are given by
\[
 T^*T^{**}=\sH \times \{0\} \quad \mbox{and} \quad T^{**}T^*=\{0\} \times \sH.
\]
\end{proposition}

Another example to illustrate what may happen with $T^*T$ and $TT^*$ if
the relation $T$ is not closed is given by the
following simple but exotic situation; this example can be seen as a
variant of Proposition~\ref{einz} in the context of relations.

\begin{example}
Let $\sR$ be a dense linear subspace of $\sK$ and define $T \in \bL(\sH, \sK)$ by
$T=\sH \times \sR$. Then  $T^*=\{0\} \times \{0\}$ and $T^{**}=\sH \times \sK$.
Now observe that
\[
T^*T=\sH \times \{0\} \quad \mbox{and} \quad TT^*=\{0\} \times \sR.
\]
Thus $T^*T$ is selfadjoint but $TT^*$ is selfadjoint if and only $\sR=\sK$,
which is equivalent to $T$ is closed. Observe for completeness that
$T^*T^{**}=\sH \times \{0\}$ and $T^{**}T^*=\{0\} \times \sK$.
\end{example}

The end of this section is devoted to the following question discussed in Research Gate in 2022:

\medskip

\textit{Is $TT^*$ (or $T^*T$) densely defined if $T$ is densely defined and symmetric?}

\medskip

\noindent
Note the requirement that $T$ be symmetric. It will be shown that the answer to this problem
in both cases is negative; see Propositions \ref{Neg1},~\ref{Neg2} below.
Block relations as in \eqref{A} can be used
to give an extension of Theorem~\ref{thm1} to the setting
where $T$ is replaced by a densely defined operator $Z$
which is essentially selfadjoint.
In fact, there exists a densely defined essentially selfadjoint operator $Z$  such that
$Z^*Z^{**}$ is a selfadjoint operator, while $Z^*Z$ is not densely defined.

\begin{proposition}\label{Neg1}
Let $T \in \bL(\sH)$ be the closable densely defined operator
from Theorem {\rm \ref{thm1}}.
Then $Z \in \bL(\sH \times \sH)$ in \eqref{A} is a densely defined essentially selfadjoint operator.
Moreover, $Z^*Z$ and $Z^*Z^{**}$ are given by
\[
Z^*Z
   =\begin{pmatrix}
     T^{**}T^* & 0 \\
     0 & \{0\} \times \{0\}
   \end{pmatrix}
\quad \mbox{and} \quad
Z^*Z^{**}
   =\begin{pmatrix}
     T^{**}T^* & 0 \\
     0 & T^*T^{**}
   \end{pmatrix},
\]
where $T^{**}T^*$ and $T^*T^{**}$ are selfadjoint operators.
In particular,  $Z^*Z$ is not densely defined while $Z^*Z^{**}$ is a selfadjoint operator.
\end{proposition}
\begin{proof}
To see this apply the formulas \eqref{A} and \eqref{Astar},
see also \eqref{A1}, to the densely defined operators
$T=B^{\half}A^{-1}$, $T^*=A^{-1}B^{\half}$,
and $T^{**}=(A^{-1}B^{\half})^*$; cf. \eqref{T1}.
\end{proof}

Notice also that
\[
 (Z^*Z)^*
   =\begin{pmatrix}
     T^{**}T^* & 0 \\
     0 & \sH \times \sH
   \end{pmatrix}.
\]
Consequently,  the selfadjoint extensions of $Z^*Z$ are given by
the componentwise orthogonal sum
\[
 Z^*Z \hoplus H,
\]
where $H$ is a selfadjoint relation in $\sH$.

\medskip

As a last result it is shown that
there exists a densely defined essentially selfadjoint operator $Z$,  which is  bounded {\rm(}and hence closable{\rm )},
such that $Z^{**}Z^*$ is a bounded selfadjoint operator, while $ZZ^{**}$ is not densely defined.

\begin{proposition}\label{Neg2}
Let $T \in \bL(\sH)$ be the bounded densely defined operator
from Proposition~{\rm \ref{thm11}}.
Then $Z \in \bL(\sH \times \sH)$ in \eqref{A} is a bounded densely defined symmetric operator.
Moreover, $ZZ^*$ and $Z^{**}Z^*$ are given by
\begin{equation*}\label{C}
ZZ^*
   =\begin{pmatrix}
     \{0\} \times \{0\} & 0 \\
     0 & T^{*}T^{**}
   \end{pmatrix}
\quad \mbox{and} \quad
 Z^{**}Z^*
   =\begin{pmatrix}
     T^{**}T^* & 0 \\
     0 & T^*T^{**}
   \end{pmatrix},
\end{equation*}
where $T^{**}T^*$ and $T^*T^{**}$ are bounded selfadjoint operators.
In particular,  $ZZ^*$ is not densely defined while $Z^{**}Z^*$ is a bounded selfadjoint operator.
\end{proposition}
\begin{proof}
To see this apply the formulas \eqref{A}, \eqref{Astar}, and \eqref{A1}
to the bounded densely defined operators operators
$T=AB^{-\half}$, $T^* =B^{-\half}A$, and $T^{**}=(B^{-\half}A)^*$
from Proposition~\ref{thm11}.
\end{proof}

\bigskip\noindent
\textbf{Acknowledgments.}
Henk de Snoo gratefully acknowledges kind invitations to
the Department of Mathematics and Statistics of the University of Vaasa, Finland.

\bigskip\noindent
\textbf{Funding Declaration.}
Seppo Hassi gratefully acknowledges support from the Research Council of Finland / Suomen Akatemia (project no. 358155).

\medskip\noindent
\textbf{Conflicts of Interest.}
The authors have no conflict of interest to declare that are relevant to the content of this article.

\medskip\noindent
\textbf{Data Availability Statement.}
No data were collected, generated or consulted in connection with this research.

 \end{document}